\documentclass[12pt,a4paper]{article}
\usepackage{amsmath}
\usepackage{amssymb}
\usepackage{amsthm}
\usepackage{float}
\usepackage{tikz}
\usepackage{multirow}
\usepackage{alltt}
\usepackage{bbm}
\usepackage{hyperref}

\tikzset{my loop/.style =  {to path={
  \pgfextra{}
  [looseness=12,min distance=10mm]
  \tikz@to@curve@path},font=\sffamily\small
  }}
\makeatletter

\theoremstyle{definition}

\theoremstyle{remark}

\begin{document}

\title{An Improved Bound for Upper Domination of Cartesian
Products of Graphs}
\author{ Yu-Yen Chien \\ Mathematics Division \\ National Center for Theoretical Science \\  \textsc{Taiwan}
\\ \texttt{\small{yychien@ncts.ntu.edu.tw}} }
\date{\today}

\maketitle

\begin{abstract}
In this paper, we prove a problem proposed by Bre\v{s}ar: for any
graphs $G$ and $H$, $\Gamma(G\square H)\ge\Gamma(G)\Gamma(H)+
\min\{|V(G)|-\Gamma(G),|V(H)|-\Gamma(H)\}$, where $\Gamma(G)$
denotes the upper domination number of $G$.
\end{abstract}

For a simple graph $G$, a subset $D$ of $V(G)$ is a {\it dominating
set} of $G$ if every vertex in $V(G)\setminus D$ has at least one
neighbor in $D$. The {\it domination number} $\gamma(G)$ is the
minimum cardinality of a dominating set of $G$. For graphs $G$ and
$H$, the {\it Cartesian product} $G\square H$ is the graph with
vertex set $V(G)\times V(H)$ and edge set
$\{(u_1,v_1)(u_2,v_2)|u_1=u_2,v_1v_2\in E(H)$, or $v_1=v_2,u_1u_2\in
E(G)\}$. Vizing \cite{a} suggested a conjecture regarding domination
in Cartesian products of graphs:

\begin{flushleft}
{\bf Vizing Conjecture.} For any graphs $G$ and $H$,
$\gamma(G\square H)\ge\gamma(G)\gamma(H)$.
\end{flushleft}

\noindent This conjecture is one of the main problems in domination
theory. See \cite{sur,e,f} for surveys and \cite{aha,g,h} for recent
progress.

For a simple graph $G$, a dominating set $D$ of $G$ is a {\it
minimal dominating set} if no proper subset of $D$ is a dominating
set of $G$. The {\it upper domination number} $\Gamma(G)$ of $G$ is
the maximum cardinality of a minimal dominating set of $G$. The
definition of the domination number $\gamma(G)$ can be rephrased as
the minimum cardinality of a minimal dominating set of $G$, and
clearly we have $\Gamma(G)\ge\gamma(G)$ for any graph $G$.
Nowakowski and Rall \cite{b} conjectured that for any graphs $G$ and
$H$, $\Gamma(G\square H)\ge\Gamma(G)\Gamma(H)$. Bre\v{s}ar \cite{c}
proved a slightly stronger bound: $\Gamma(G\square
H)\ge\Gamma(G)\Gamma(H)+1$ for any nontrivial graphs $G$ and $H$,
where a {\it nontrivial} graph is a graph with at least one edge.
Bre\v{s}ar also proposed the following question: does the inequality
$$\Gamma(G\square H)\ge\Gamma(G)\Gamma(H)+\min\{|V(G)|-\Gamma(G),|V(H)|-\Gamma(H)\}$$
hold for any graphs $G$ and $H$? We prove this inequality in this
paper.

We examine some basic properties about minimal dominating set first.
$G$ is a graph and $D$ is a dominating set of $G$. We say a vertex
in $D$ is $D$-{\it isolated} if it is not adjacent to any other
vertices in $D$. We say a vertex $v \in D$ has a {\it private}
$D$-{\it neighbor} $u$ if $u\notin D$ is a neighbor of $v$ and is no
neighbor of any other vertices in $D$. We have the following
fundamental results.

\begin{flushleft}
{\bf Lemma 1.}\cite{d} A dominating set $D$ of a graph $G$ is a
minimal dominating set if and only if every vertex in $D$ is
$D$-isolated or has a private $D$-neighbor.
\end{flushleft}

\begin{flushleft}
{\bf Lemma 2.} Given a dominating set $D$ of a graph $G$, we can
always find a subset $D'$ of $D$ such that $D'$ is a minimal
dominating set of $G$. Moreover, if $v\in D$ is $D$-isolated or has
a private $D$-neighbor, then $v$ must be in $D'$.
\end{flushleft}

For the graph $G\square H$, $u\in V(G)$, $v\in V(H)$ , we call
$\{u\}\times V(H)$ {\it row} $u$ and $V(G)\times \{v\}$ {\it column}
$v$. In the following, we consider upper domination in the Cartesian
products of some basic graphs.

\begin{flushleft}
{\bf Proposition 3.} For $m,n\ge 2$, $\Gamma(K_m\square K_n)=$
max$\{m,n\}$.
\end{flushleft}
\noindent {\bf Proof.} Assume $n\ge m$. It is clear that a row of
$K_m\square K_n$ is a minimal dominating set, and $\Gamma(K_m\square
K_n)\ge n$. Suppose $K_m\square K_n$ has a minimal dominating set
$D$ with cardinality greater than $n$. If every column of
$K_m\square K_n$ has a vertex in $D$, the $D=n$. Therefore, there
exists a column without a vertex in $D$. To dominate vertices in
this column, every row of $K_m\square K_n$ must have a vertex in
$D$. Then, we have $D=m$, a contradiction. \qed

The following proposition is implicitly proved by Gutin and
Zverovich \cite{i}.

\begin{flushleft}
{\bf Proposition 4.} For any graph $G$, $\Gamma(K_2\square G)=|G|$.
\end{flushleft}
\noindent {\bf Proof.} We prove the case that $G$ is connected and
$|V(G)|\ge 2$, and then the general case follows. Since a row of
$K_2\square G$ is a minimal dominating set, we have
$\Gamma(K_2\square G)\ge |G|$. Suppose $D$ is a minimal dominating
set of $K_2\square G$. For $v\in D$ having a private $D$-neighbor,
assign one of its private $D$-neighbor as $v'$. For $v\in D$ having
no private $D$-neighbor, assign the other vertex in column $v$ as
$v'$. For any $v\in D$, $v'\notin D$. For any different $u$, $v$ in
$D$, $u'\neq v'$. Therefore, $|D|\leq\frac{|K_2\square G|}{2}=|G|$,
and we have  $\Gamma(K_2\square G)=|G|$. \qed

\begin{flushleft}
{\bf Proposition 5.} For $m\ge 3$ and $n\ge 2$, $\Gamma(K_m\square
K_{1,n})=m+n-2$.
\end{flushleft}
\noindent {\bf Proof.} Let $V(K_m)=\{u_1,...,u_m\}$ and
$V(K_{1,n})=\{v_0,v_1,...,v_n\}$, where $v_0$ is the vertex of
degree $n$. $\{(u_i,v_1),(u_1,v_j)|m\ge i\ge 2,n\ge j\ge 2\}$ is a
minimal dominating set of $K_m\square K_{1,n}$, and we have
$\Gamma(K_m\square K_{1,n})\ge m+n-2$. Suppose $K_m\square K_{1,n}$
has a minimal dominating set $D$ with $|D|>m+n-2$. If $D$ has a
vertex in each column of $K_m\square K_{1,n}$, then $|D|=n+1>m+n-2$,
violating the assumption that $m\ge 3$. Therefore, there exists
$v_j$ such that column $v_j$ contains no vertex in $D$. If $j\neq
0$, $(u_i,v_0)$ must be in $D$ to dominate $(u_i,v_j)$, and we have
$D=m>m+n-2$, violating the assumption that $n\ge 2$. Therefore,
column $v_0$ contains no vertex in $D$. To dominate $(u_i,v_0)$,
there exists some $(u_i,v_j)\in D$. Assign one of these vertices as
$(u_i,v_0)'$. Let $D'=\{(u_i,v_0)'|m\ge i\ge 1\}$. Notice that $D'$
can not be a column. Columns containing a vertex in $D'$ do not
contain any vertex in $D\setminus D'$. Except column $v_0$, columns
containing no vertex in $D'$ contains exact one vertex in
$D\setminus D'$, and there are at most $n-2$ such columns.
Therefore, $|D|\le m+n-2$, a contradiction. \qed

\begin{flushleft}
{\bf Proposition 6.} For $l,m,n\ge 2$, we have $$\Gamma(K_l\square
K_{m,n})=\max\{m+n,2l,m+l-2,n+l-2\}.$$
\end{flushleft}
\noindent {\bf Proof.} Suppose that $V(K_l)=\{u_1,...,u_l\}$, and
$K_{m,n}$ has partite sets $V=\{v_1,...,v_m\}$ and $V'=
\{v_1',...,v_n'\}$. Row $u_1$, union of column $v_1$ and column
$v_1'$, $\{(u_i,v_1),(u_1,v_j)|l\ge i\ge 2, m\ge j\ge 2\}$, and
$\{(u_i,v_1'),(u_1,v_j')|l\ge i\ge 2, n\ge j\ge 2\}$ are minimal
dominating sets of $K_l\square K_{m,n}$, so we have
$\Gamma(K_l\square K_{m,n})\ge$ max$\{m+n,2l,m+l-2,n+l-2\}$. Suppose
$D$ is a minimal dominating set of $K_l\square K_{m,n}$. If each
column of $K_l\square K_{m,n}$ has a vertex in $D$, then $|D|=m+n$.
If each of $\{u_i\}\times V$ and $\{u_i\}\times V'$ has a vertex in
$D$, then $|D|=2l$. Suppose $|D|\ge$ max$\{m+n,2l\}$, and then there
exists a column $v^*$ without vertex in $D$. If $v^*\in V'$, each of
$\{u_i\}\times V$ must have a vertex $w_i$ in $D$, and there exists
some $\{u_i\}\times V'$ without vertex in $D$. Without loss of
generality, suppose $\{u_1\}\times V'$ has no vertex in $D$ and
$w_1=(u_1,v_1)$. To dominate $(u_1,v_j)$ where column $v_j$ has no
vertex in $\{w_i\}$, column $v_j$ must have a vertex $x_j\in D$.
Notice that $\{w_i\}\cup\{x_j\}$ is a dominating set of $K_l\square
K_{m,n}$. If $\{w_i\}$ is a column, then $\{w_i\}\cup\{x_j\}$ can
not be minimal. So $\{w_i\}$ is not a column, $|\{x_j\}|\le m-2$,
and $|D|\le m+l-2$. If $v^*\in V$, similarly we have $|D|\le n+l-2$.
\qed

We define two graphs $X_n$ and $X_n'$ here. $X_n$ has vertex set
$\{u_0,u_1,...,u_n\}\cup\{v_0,v_1,...,v_n\}$ and edge set
$\{u_iu_j|i>j\}\cup\{v_iv_j|i>j\}\cup\{u_iv_i|i\neq 0\}$. We call
the subgraph of $X_n$ induced by $\{u_0,u_1,...,u_n\}$ an {\it upper
cell}, the subgraph induced by $\{v_0,v_1,...,v_n\}$ a {\it lower
cell}. We use $X_n'$ to denote the subgraph induced by
$\{u_1,...,u_n\}\cup\{v_0,v_1,...,v_n\}$ in $X_n$. Consider
$\Gamma(X_n)$ and $\Gamma(X_n')$. For $X_n$, every cell of $X_n$
must have a vertex to dominate $u_0$ and $v_0$, and these two
vertices dominate $X_n$. Therefore, we have $\Gamma(X_n)=2$. For
$X_n'$, suppose $D$ is a minimal dominating set of $\Gamma(X_n')$.
If each of $\{u_1,...,u_n\}$ and $\{v_0,v_1...,v_n\}$ has a vertex
in $D$, then $|D|=2$. If $D>2$, $\{v_0,v_1...,v_n\}$ must have a
vertex in $D$ to dominate $v_0$ and therefore $\{u_1,...,u_n\}$ has
no vertex in $D$. To dominate $u_i$, $v_i$ must be in $D$. We have
$D=\{v_1...,v_n\}$ and $\Gamma(X_n')=n$. $X_n'$ is the induced
subgraph of $X_n$ and $\Gamma(X_n')$ is greater than $\Gamma(X_n)$.
Similar situation also happens in Cartesian products of graphs. The
following proposition gives an example.

\begin{flushleft}
{\bf Proposition 7.} For sufficiently large $n$, we have
$\Gamma(X_n'\square K_3)>\Gamma(X_n\square K_3)$.
\end{flushleft}
\noindent {\bf Proof.} Since $\{v_1,...,v_n\}\times V(K_3)$ is a
minimal dominating set of $X_n'\square K_3$, we have
$\Gamma(X_n'\square K_3)\ge 3n$. Suppose $D$ is a minimal dominating
set of $X_n\times K_3$ with $|D|\ge 3n$. If a cell in $X_n\times
K_3$ has two or more vertices in $D$, then these vertices must have
private $D$-neighbors, which must be in cells without vertex in $D$.
If $n$ is sufficiently large, it is easy to see that there are three
cells without vertex in $D$ and almost all vertices in the other
three cells are in $D$. Cells without vertex in $D$ are not all
upper cells, or $\{u_0\}\times V(K_3)$ can not be dominated.
Similarly, cells without vertex in $D$ are not all lower cells.
Suppose two upper cells have no vertex in $D$. In the remaining
upper cell, a vertex in $D$ dominates at least two vertices in cells
without vertex in $D$ and prohibits these vertices from being
private $D$-neighbors of other vertices in $D$. It follows that
$|D|<3n$, a contradiction. Therefore, we have $\Gamma(X_n'\square
K_3)\ge3n>\Gamma(X_n\square K_3)$. \qed

Suppose $x$ is a vertex having only one neighbor $y$ in a graph $G$.
If $D$ is a minimal dominating set of $G$, $D$ must contain exactly
one vertex of $x$ and $y$. To minimize $D$, it is a better choice to
put $y$ in $D$. However, if we want to maximize $D$, putting $x$ in
$D$ is not always better. For example, let $G$ be the graph with
vertex set $\{x,y\}\cup V(X_n)$ and edge set $\{xy,yu_0\}\cup
E(X_n)$, and $D$ is a minimal dominating set of $G$. Using similar
argument in determining $\Gamma(X_n)$ and $\Gamma(X_n')$, we know
$|D|=3$ if $D$ contains $x$, and $|D|=n+1$ if $D$ contains $y$.

Given a graph $G$, if we add some edges in $G$ to get a graph $G'$,
we always have $\gamma(G)\ge\gamma(G')$. However, upper domination
in graphs does not have this property. For example, let $G'$ be the
graph obtained from adding edge $u_0v_0$ in $G=X_n$. Then $G'$ is
isomorphic to $K_2\square K_{n+1}$ and $\Gamma(G)=2<\Gamma(G')=n+1$
for $n>1$. The above examples show that upper domination in graphs
probably does not behave as well as one may expect, and the argument
of minimal counterexample probably does not work well in solving
upper domination problems.

It is easy to see that a maximal independent set is a minimal
dominating set, and therefore we have $\Gamma(G)\ge\alpha(G)$ for
any graph $G$. It is clear that $\alpha(X_n)=\alpha(X_n')=2$. Since
$\Gamma(X_n)=2$ and $\Gamma(X_n')=n$, we know the inequality
$\Gamma(G)\ge\alpha(G)$ is sharp, and the difference between
$\Gamma(G)$ and $\alpha(G)$ could be quite large.

\begin{flushleft}
{\bf Lemma 8.} $G$ and $H$ are two arbitrary graphs. $I_G$ and $I_H$
are the maximal independent sets of $G$ and $H$ respectively. Let
$G'=G[V(G)\setminus I_G]$ and $H'=G[V(H)\setminus I_H]$. We have
$\Gamma(G\square H)\ge |I_G||I_H|+\Gamma(G'\square H')$.
\end{flushleft}
\noindent {\bf Proof.} Suppose $D_1=I_G\times I_H$ and $D_2$ is the
maximum minimal dominating set of $G'\square H'$. We claim that
$D_1\cup D_2$ is a minimal dominating set of $G\square H$ and then
we have the inequality. Since $I_G$ and $I_H$ dominate $V(G)$ and
$V(H)$ respectively, $D_1$ dominates $I(G)\times V(H)$ and
$V(G)\times I(H)$, and $D_1\cup D_2$ dominates $V(G)\times V(H)$.
Notice that there is no edge between $D_1$ and $D_2$. Vertices in
$D_1$ are $(D_1\cup D_2)$-isolated, and $D_2$-isolated vertices are
still $(D_1\cup D_2)$-isolated. If a vertex $u$ in $D_2$ has a
private $D_2$-neighbor $v$ in $G'\square H'$, $v$ is not adjacent to
any vertex in $D_1$ and is still a private $(D_1\cup D_2)$-neighbor
of $u$ in $G\square H$. By Lemma 1, we know that $D_1\cup D_2$ is a
minimal dominating set of $G\square H$. \qed

In $G\square H$, we can choose vertices from different rows and
columns to form an independent set. Therefore, we have
$\Gamma(G\square H)\ge\alpha(G\square H)\ge $
min$\{|V(G)|,|V(H)|\}$. The following lemma improves this result.

\begin{flushleft}
{\bf Lemma 9.} For any nontrivial graph $G$ and arbitrary graph $H$,
we have $\Gamma(G\square H)\ge|V(H)|$, and the equality holds only
if $G$ is a complete graph or $K_{1,2}$.
\end{flushleft}
\noindent {\bf Proof.} Since $G$ is nontrivial, we can always find a
component $G_0$ of $G$ with $|V(G_0)|\ge 2$. Choose a vertex $u$ in
$V(G_0)$. The set $N(u)$ of neighbors of $u$ in $G_0$ is not empty.
In case row $u$ is a dominating set of $G_0\square H$, it is a
minimal dominating set of $G_0\square H$ and we have
$\Gamma(G\square H)\ge\Gamma(G_0\square H)\ge|V(H)|$. Suppose row
$u$ is not a dominating set of $G_0\square H$. Let
$D=(V(G_0)\setminus N(u))\times V(H)$. Clearly, $D$ is a dominating
set of $G_0\square H$, and we can find a subset $D'$ of $D$ such
that $D'$ is a minimal dominating set of $G_0\square H$. For every
vertex $v\in V(H)$, $D'$ must have a vertex in column $v$ of
$G_0\square H$ to dominate vertices in $N(u)\times \{v\}$, and we
have $\Gamma(G\square H)\ge\Gamma(G_0\square H)\ge|V(H)|$.

Suppose $\Gamma(G\square H)=|V(H)|$. $G$ must be connected.
Otherwise, we have $\Gamma(G\square H)>\Gamma(G_0\square
H)\ge|V(H)|$. If $G$ is not bipartite, an induced subgraph of $G$
from removing a maximal independent set of $G$ is nontrivial. If $G$
is bipartite but not complete bipartite, choose two vertices $u$ and
$u'$ from different partite sets of $G$ such that $u$ and $u'$ are
not adjacent. Suppose $I$ is a maximal independent set of $G$
containing $u$ and $u'$. Then $G[V(G)\setminus I]$ must be
nontrivial. Otherwise, $G$ is not connected. Therefore, if $G$ is
not a complete bipartite graph, we can always find a maximal
independent set $I_G$ of $G$ such that $G'=G[V(G)\setminus I_G]$ is
nontrivial. Suppose $I_H$ is a maximal independent set of $H$ and
$H'=H[V(H)\setminus I_H]$. By Lemma 8 and the inequality proven
above, we have $\Gamma(G\square H)=|V(H)|\ge
|I_G||I_H|+\Gamma(G'\square H')\ge |I_G||I_H|+|V(H')|$. It implies
$|I_G|=1$, which means $G$ is a complete graph. Proposition 3 shows
that the equality does hold while $H$ is also a complete graph with
$V(H)\ge V(G)$.

Consider the case that $G$ is complete bipartite. If $G=K_{1,1}$,
which is also a complete graph, Proposition 4 shows that the
equality always holds. If $G=K_{1,2}$, Proposition 5 shows that the
equality does hold while $H$ is a complete graph with $V(H)\ge 3$.

If $G=K_{1,m}$, $m\ge 3$, let $V(G)=\{u_0,u_1,...,u_m\}$, where
$u_0$ is the vertex of degree $m$. If $H$ has a vertex of degree
$0$, it is easy to see that $\Gamma(G\square H)>|V(H)|$. Suppose
every vertex of $H$ has degree at least 1, and $D_H$ is a minimal
dominating set of $H$. Let $R_H=V(H)\setminus D_H$. Notice that
$R_H$ is not empty, and every vertex in $D_H$ has a neighbor in
$R_H$. In other words, $R_H$ is a dominating set of $H$. By Lemma 1,
it is easy to check that $(\{u_m\}\times
R_H)\cup(\{u_1,...,u_{m-1}\}\times D_H)$ is a minimal dominating set
of $G\square H$, and we have $\Gamma(G\square H)>|V(H)|$. If
$G=K_{m,n}$, $m,n\ge2$, let $u$ and $u'$ be two vertices in
different partite sets of $G$. Then $\{u,u'\}\times V(H)$ is a
minimal dominating set of $G\square H$, and we have $\Gamma(G\square
H)>|V(H)|$. \qed

Let $D$ be the minimal dominating set of a graph $G$. We define
$D^P$ as the set of vertices in $D$ having a private $D$-neighbor,
and $D^I=D\setminus D^P$. By Lemma 1, we know vertices in $D^I$ are
$D$-isolated. Now we are ready to prove the main theorem.

\begin{flushleft}
{\bf Theorem 10.} For any graphs $G$ and $H$, we have
$$\Gamma(G\square H)\ge\Gamma(G)\Gamma(H)+
\min\{|V(G)|-\Gamma(G),|V(H)|-\Gamma(H)\}.$$
\end{flushleft}
\noindent {\bf Proof.} If $G$ or $H$ has no edge, then
$|V(G)|=\Gamma(G)$ or $|V(H)|=\Gamma(H)$, and we have
$\Gamma(G\square H)=\Gamma(G)\Gamma(H)=\Gamma(G)\Gamma(H)+$
min$\{|V(G)|-\Gamma(G),|V(H)|-\Gamma(H)\}$. Now we prove the case
that every vertex of $G$ and $H$ has degree at least 1, and then the
case that $G$ and $H$ are nontrivial follows. Let $D_G$ and $D_H$ be
the minimal dominating sets of $G$ and $H$ respectively. Define
$R_G=V(G)\setminus D_G$, $R_H=V(H)\setminus D_H$, $G'=G[R_G]$,
$H'=H[R_H]$. If $D_G^P$ and $D_H^P$ are both empty, then $D_G$ and
$D_H$ are maximal independent sets of $G$ and $H$ respectively. By
Lemma 8, we have $\Gamma(G\square H)\ge|D_G||D_H|+\Gamma(G'\square
H')\ge|D_G||D_H|+$ min$\{|V(G)|-|D_G|,|V(H)|-|D_H|\}$.

Suppose $D_G^P$ is not empty. For each vertex $(u,v)\in D_G^P\times
D_H$, we select a corresponding vertex $(u,v)'=(u',v)$ such that
$u'$ is a private $D_G$-neighbor of $u$ in $G$. Let
$D_1=\{(u,v)'|(u,v)\in D_G^P\times D_H\}$. For each vertex $(u,v)\in
D_G^I\times D_H^P$, we select a corresponding vertex
$(u,v)''=(u,v'')$ such that $v''$ is a private $D_H$-neighbor of $v$
in $H$. Let $D_2=\{(u,v)''|(u,v)\in D_G^I\times D_H^P\}$. Since
every vertex of $G$ has degree at least 1, every vertex in $D_G^I$
must have a neighbor in $R_G$, and therefore $R_G$ is a dominating
set of $G$. Similarly, $R_H$ is a dominating set of $H$. It is easy
to check that $D=D_1\cup D_2\cup(D_G^I\times D_H^I)\cup(R_G\times
R_H)$ is a dominating set of $G\square H$. Let $D'\subset D$ be a
minimal dominating set of $G\square H$. Notice that each vertex in
$D_G^I\times D_H^I$ is $D$-isolated, and each vertex in $D_1\cup
D_2$ has a private $D$-neighbor in $G\square H$. By Lemma 2, we have
$D_1\cup D_2\cup(D_G^I\times D_H^I)\subset D'$. There is no edge
between $D_1\cup D_2\cup(D_G^I\times D_H^I)$ and $D_G^P\times R_H$.
For each $v\in R_H$, $D'$ must have at least one vertex in
$R_G\times \{v\}$ to dominate vertices in $D_G^P\times \{v\}$.
Therefore, we have $\Gamma(G\square
H)\ge|D'|\ge|D_1|+|D_2|+|D_G^I\times
D_H^I|+|R_H|=|D_G||D_H|+|V(H)|-|D_H|$. It is easy to check that the
above argument works no matter $D_G^I$, $D_H^I$, $D_H^P$ are empty.
\qed

In the above proof, if $D^P_G$ and $D^P_H$ are both empty, then we
have $\Gamma(G\square H)\ge|D_G||D_H|+\Gamma(G'\square H')$. In case
that $G'$ and $H'$ are nontrivial, such as $G$ and $H$ are not
bipartite, then we have $\Gamma(G\square H)\ge|D_G||D_H|+$
max$\{|V(G)|-|D_G|,|V(H)|-|D_H|\}$ by Lemma 9. If $D^P_G$ is not
empty, we have $\Gamma(G\square H)\ge|D_G||D_H|+|V(H)|-|D_H|$.
Therefore, if $D^P_G$ and $D^P_H$ are not empty, then we also have
$\Gamma(G\square H)\ge|D_G||D_H|+$
max$\{|V(G)|-|D_G|,|V(H)|-|D_H|\}$. The above discussion shows that
the inequality in Theorem 10 can be improved to $\Gamma(G\square
H)\ge\Gamma(G)\Gamma(H)+$ max$\{|V(G)|-\Gamma(G),|V(H)|-\Gamma(H)\}$
in many cases. However, in general the minimum in the inequality can
not be replaced by maximum. Let $G=K_6$ and $H=K_{4,8}$. We know
$\Gamma(G)=1$, $\Gamma(H)=8$, and $\Gamma(G\square H )=12$ by
Proposition 6. It shows that the inequality in Theorem 10 is sharp.
We wonder if the lower bound of $\Gamma(G\square H)$ is close to
$\Gamma(G)\Gamma(H)+$ max$\{|V(G)|-\Gamma(G),|V(H)|-\Gamma(H)\}$.
For example, for $G=K_l$ and $H=K_{m,n}$, $l,m,n\ge 2$, we have
$\Gamma(G\square H)\ge\Gamma(G)\Gamma(H)+$
max$\{|V(G)|-\Gamma(G),|V(H)|-\Gamma(H)\}-1$ by Proposition 6.

We give an alternative lower bound of $\Gamma(G\square H)$ to end
this paper. This lower bound can be much better than Theorem 10 in
some cases. For a graph $G$ and $S\subset V(G)$, we define
$N[S]=S\cup \{v\in V(G)| v$ has a neighbor in $S\}$.

\begin{flushleft}
{\bf Theorem 11.} $G$ and $H$ are two arbitrary graphs, and $D$ is a
minimal dominating set of $H$. We have $\Gamma(G\square H)\ge
|V(G)||D^P|+\gamma(G)|D^I|$. Moreover, if $N[D^P]\cup D^I=V(H)$,
then we have $\Gamma(G\square H)\ge |V(G)||D^P|+\Gamma(G)|D^I|$.
\end{flushleft}
\noindent {\bf Proof.} Let $D_1=V(G)\times D^P$ and $D_2=V(G)\times
D^I$. Clearly, $D_1\cup D_2$ is a dominating set of $G\square H$,
and suppose $D_3\subset D_1\cup D_2$ is a minimal dominating set of
$G\square H$. Notice that vertices in $D_1$ have private $(D_1\cup
D_2)$-neighbors in $G\square H$, and therefore $D_1\subset D_3$. For
$v\in D^I$, $D_3$ must have at least $\gamma(G)$ vertices in column
$v$ to dominate column $v$. Therefore, $\Gamma(G\square
H)\ge|D_3|\ge|V(G)||D^P|+\gamma(G)|D^I|$. In case $N[D^P]\cup
D^I=V(H)$, let $D'$ be a maximum minimal dominating set of $V(G)$.
It is easy to check that $D_1\cup(D'\times D^I)$ is a minimal
dominating set of $G\square H$, and we have $\Gamma(G\square H)\ge
|V(G)||D^P|+\Gamma(G)|D^I|$. \qed

\vspace{2mm}
\begin{flushleft}
{\bf Acknowledgements.} The author would like to thank the support
of National Center for Theoretical Science, Taiwan.
\end{flushleft}

\end{document}